\documentclass[reqno,12pt]{amsart}
\usepackage{epsfig}
\usepackage{amsfonts,amsmath,amssymb,mathrsfs,verbatim,amsthm,amscd}
\def\II{\hbox{{1}\kern-.25em\hbox{l}}}
\newcommand{\V}{\mathbb V}
\newcommand{\be}{\begin{equation}}
\newcommand{\ee}{\end{equation}}
\newcommand{\ba}{\begin{eqnarray}}
\newcommand{\ea}{\end{eqnarray}}

\newcommand{\C}{\mathbb C}

\newcommand{\Z}{\mathbb Z}

\newcommand{\T}{\mathbb T}

\newcommand{\eg}{\Gamma}
\newtheorem{theorem}{Theorem}

\begin{document}

\title[Aspects of elliptic hypergeometric functions]
{Aspects of elliptic hypergeometric functions}

\author{V.P. Spiridonov}
 \address{JINR, Dubna, Moscow region, Russia}

\thanks{Talk given at the International Conference
``The Legacy of Srinivasa Ramanujan", 17--22 December 2012, Delhi, India}

\begin{abstract}
General elliptic hypergeometric
functions are defined by elliptic hypergeometric integrals.
They comprise the elliptic beta integral, elliptic analogues of the
Euler-Gauss hypergeometric function and Selberg
integral, as well as elliptic extensions of many other
plain hypergeometric and $q$-hypergeometric constructions.
In particular, the Bailey chain technique, used for
proving Rogers-Ramanujan type identities, has been
generalized to integrals. At the elliptic level it
yields a solution of the Yang-Baxter equation
as an integral operator with an elliptic hypergeometric
kernel. We give a brief survey of the developments in this field.

\end{abstract}

\maketitle

\tableofcontents

\section{Elliptic hypergeometric integrals}

Hypergeometric functions lie at the center of the world of special functions \cite{aar}.
Ramanujan obtained many important results in the theory of hypergeometric functions
and their $q$-analogues. It is therefore natural to give at his
jubilee conference a survey of the top known special functions of
hypergeometric type --
the elliptic hypergeometric functions. General representatives of
these functions are defined by elliptic hypergeometric integrals
introduced by the author in 2000 \cite{spi:umn} and a general
setup of their theory was formulated in \cite{spi:theta2,spi:thesis}.
An overview of the results obtained prior to 2008 is given in \cite{spi:umnrev}.

In the univariate case the contour integrals
$$
I=\int_C\Delta(u)du
$$
are called elliptic hypergeometric integrals, if the (meromorphic) kernel function $\Delta(u)$
satisfies a first order finite difference equation
$$
\Delta(u+\omega_1)=h(u;\omega_2,\omega_3)\Delta(u),
$$
where $h(u;\omega_2,\omega_3)$ is an elliptic function,
$$
h(u+\omega_2)=h(u+\omega_3)=h(u),\quad \text{Im}(\omega_2/\omega_3)\neq 0,
$$
and $\omega_{1,2,3}$ are some (in general incommensurate) complex variables.
Define two bases
$$
p  =  e^{2 \pi \textup{i}
\omega_3/\omega_2}, \quad q  =  e^{2 \pi \textup{i} \omega_1/\omega_2}
$$
and demand that $\Delta(u):=\rho(z)$ is a meromorphic function
in the variable
$z=e^{2\pi\textup{i} u/\omega_2}\in\C^*$. Then we can write
$$
I=\int\rho(z)\frac{dz}{z}, \qquad \rho(qz)=h(z;p)\rho(z),\quad h(pz;p)=h(z;p),
$$
where
$$
h(z;p)=\prod_{k=1}^m\frac{\theta(t_kz;p)}{\theta(w_kz;p)},
\qquad \prod_{k=1}^mt_k=\prod_{k=1}^mw_k,
$$
is an arbitrary elliptic function defined as a ratio of products of
theta functions
$$
\theta(z;p)=(z;p)_\infty(pz^{-1};p)_\infty,
\qquad (z;p)_\infty=\prod_{j=0}^\infty(1-zp^j), \quad |p|<1.
$$
The integer parameter $m\geq 2$ is called the order of the elliptic function
$h(z;p)$ and $t_k, w_k$ are arbitrary parameters fixing its divisor.

Now, due to the factorization of $h(z;p)$, it is sufficient to solve the
following linear first order $q$-difference equation
$$
f(qz)=\theta(z;p)f(z).
$$
Its particular solution is given by the (standard) elliptic gamma function
\begin{equation}
\Gamma(z;p,q)
=\prod_{j,k=0}^\infty\frac{1-z^{-1}p^{j+1}q^{k+1}}{1-zp^{j}q^{k}},\qquad |p|, |q|<1.
\label{eg}\end{equation}

The multiple gamma functions were proposed by Barnes \cite{bar:multiple}.
Known generalizations of Euler's gamma function, including \eqref{eg},
can be built as combinations of Barnes' gamma functions.
Various properties of the elliptic gamma functions were investigated by
Jackson \cite{jac:basic}, Baxter \cite{bax}, Ruijsenaars (who coined its name)
\cite{rui:first}, Felder and Varchenko \cite{fel-var:elliptic},
the author \cite{spi:theta2}, and Rains \cite{rai:limits}.

As a result, for $|q|, |p|<1$ we find the general form of elliptic
hypergeometric integrals as \cite{spi:theta2,spi:thesis}
$$
I=\int\prod_{k=1}^m\frac{\Gamma(t_kz;p,q)}{\Gamma(w_kz;p,q)}\frac{dz}{z},
\qquad \prod_{k=1}^mt_k=\prod_{k=1}^mw_k.
$$

In the following we use the conventions
\begin{eqnarray*}
&& \Gamma(t_1,\ldots,t_n;p,q):=\Gamma(t_1;p,q)\cdots\Gamma(t_n;p,q),\quad
\\ &&
\Gamma(tz^{\pm k};p,q):=\Gamma(tz^k;p,q)\Gamma(tz^{-k};p,q),
\\ &&
\Gamma(tx^{\pm1}z^{\pm1};p,q):=\Gamma(t x z^{\pm1};p,q)\Gamma(t x^{-1}z^{\pm1};p,q).
\end{eqnarray*}

\section{The elliptic beta integral}

The key theorem for elliptic hypergeometric integrals was proved in \cite{spi:umn}.

\begin{theorem} Let $|p|, |q|, |t_j|<1$, $\prod_{j=1}^6t_j=pq$. Then
\begin{equation}
\frac{(p;p)_\infty(q;q)_\infty}{4\pi \textup{i}}
\int_\T\frac{\prod_{j=1}^6
\eg(t_jz^{{\pm 1}};p,q)}{\eg(z^{\pm 2};p,q)}\frac{dz}{z}
=\prod_{1\leq j<k\leq6}\eg(t_jt_k;p,q),
\label{ellbeta}\end{equation}
where $\T$ is the unit circle with positive orientation.
\end{theorem}

This was a fundamentally new exactly computable integral with the
following properties.

\begin{itemize}

\item
It represents the top known hypergeometric generalization
of Newton's binomial theorem and its $q$-analogue \cite{aar},
i.e. it can be called the elliptic binomial theorem.

\item It desribes the top known generalization of Euler's beta integral
 \cite{aar}, including at the intermediate steps the Askey-Wilson \cite{aw}
and Rahman  \cite{rah:integral} $q$-beta integrals.

\item
Formula \eqref{ellbeta} obeys the $W(E_6)$ group of symmetries -- the Weyl
group of the exceptional root system $E_6$  (see \cite{SV1}
for a detailed discussion of this property).

\item
This integral defines the orthogonality measure for two-index elliptic
biorthogonal functions \cite{spi:theta2}. The latter functions
represent an elliptic extension of the famous Askey-Wilson
orthogonal polynomials \cite{aw} and Rahman biorthogonal rational
functions \cite{rah:integral}. They constitute also the
continuous measure extension of the discrete elliptic biorthogonal
rational functions of Zhedanov and the author \cite{spi-zhe:spectral}
(elliptic analogues of Wilson's functions \cite{wil:orthogonal}).

\item
In a special limit this exact integral evaluation formula can be degenerated
to the Frenkel-Turaev summation formula (an elliptic functions identity of
hypergeometric form) \cite{ft}. The latter sum can be degenerated to
the terminating Jackson $_8\varphi_7$-sum and further on to the Dougall
$_7F_6$-sum, which was discovered also by Ramanujan.  In general,
elliptic hypergeometric series emerge as residue sums for particular
sequences of poles of elliptic hypergeometric integral kernels.
Such series play an important role \cite{spi:umnrev}, but we do not describe
them for brevity.

\item
At the univariate level only one formula of such type has been found so far,
but there are many multidimensional extensions to computable
integrals on root systems
\cite{bult:trafo,die-spi:elliptic,die-spi:selberg,rai:trans,rai:lit,spi:theta2,
spi:short,SV2,SV3,spi-war:inversions}.

\end{itemize}

The first proof of this theorem used an elliptic extension of Askey's method
described in the proceedings of a Ramanujan centennial meeting \cite{ask:beta}.
Its shortest known proof is described in \cite{spi:short}.

\section{An elliptic analogue of the Euler-Gauss hypergeometric function}

The following elliptic extension of the Euler-Gauss hypergeometric function
$_2F_1(a,b;c;x)$ \cite{aar} was introduced in \cite{spi:theta2}
\begin{equation}
V(t_1,\dots,t_8;p,q)=\frac{(p;p)_\infty(q;q)_\infty}{4\pi \textup{i}}
\int_\T\frac{\prod_{j=1}^8\Gamma(t_jx^{\pm 1};p,q)}
{\Gamma(x^{\pm2};p,q)}\frac{dx}{x}, \quad \prod_{j=1}^8t_j=(pq)^2.
\label{V}\end{equation}
Here it is assumed that $|t_j|<1$, but by an appropriate change of the integration
contour this function can be analytically continued to $t_j\in\C^*$.

In contrast to the elliptic beta integral \eqref{ellbeta}, this $V$-function
obeys the $W(E_7)$-group of symmetries with the key non-trivial transformation law
discovered in  \cite{spi:theta2}
\begin{equation}
V(t;p,q)=\prod_{1\le j<k\le 4}\Gamma(t_jt_k,t_{j+4}t_{k+4};p,q)\,
V(s;p,q),
\label{E7}\end{equation}
where $|t_j|, |s_j|<1$ and
$$
\left\{
\begin{array}{cl}
s_j =\varepsilon t_j,&   j=1,2,3,4  \\
s_j = \varepsilon^{-1} t_j, &    j=5,6,7,8
\end{array}
\right.;
\quad \varepsilon=\sqrt{\frac{pq}{t_1t_2t_3t_4}}
=\sqrt{\frac{t_5t_6t_7t_8}{pq}}.
$$

An elliptic analogue of the hypergeometric equation \cite{aar} has the form
\cite{spi:thesis,spi:cs}
\begin{eqnarray}\nonumber
&& \makebox[-2em]{}
U( t ;q,p)+\mathcal{A}(t_1,t_2,\ldots,t_8,q;p)\Big(U(qt_1,q^{-1}t_2;p,q)-U( t ;p,q)\Big)
\\ &&
+\mathcal{A}(t_2,t_1,\ldots,t_8,q;p)\Big(U(q^{-1}t_1,qt_2,;p,q)-U( t ;p,q)\Big)
=0,
\label{EHE}\end{eqnarray}
where we indicate particular scaled parameters $q^{\pm1}t_j$ in the set
$t=(t_1,\ldots,t_8)$ and
$$
 \mathcal{A}(t_1,\ldots, t_8,q;p):=\frac{\theta(t_1/qt_3,t_3t_1,t_3/t_1;p)}
                 {\theta(t_1/t_2,t_2/qt_1,t_1t_2/q;p)}
\prod_{k=4}^8\frac{\theta(t_2t_k/q;p)}{\theta(t_3t_k;p)},
$$
$$
U( t ; p,q):=\frac{V( t ;p,q)}
{\prod_{k=1}^2 \Gamma(t_kt_3^{\pm 1};p,q)}.
$$
A detailed consideration of the limiting functions obtained from
the $V$-function by various degenerations is given in \cite{brs}.

\section{The elliptic Selberg integral}

The following exact integration formula was suggested by van Diejen
and the author in \cite{die-spi:elliptic} (for $n=1$ it reduces to \eqref{ellbeta}):
\begin{eqnarray}\nonumber &&
\frac{(p;p)_\infty^n(q;q)_\infty^n}{2^n n!(2\pi \textup{i})^n}
\int_{\T^n} \prod_{1\leq j<k\leq n}
\frac{\eg(tz_j^{\pm 1} z_k^{\pm 1};p,q)}{\eg(z_j^{\pm 1} z_k^{\pm 1};p,q)}
\prod_{j=1}^n\frac{\prod_{m=1}^6\eg(t_mz_j^{\pm 1};p,q)}{\eg(z_j^{\pm2};p,q)}
\frac{dz_1}{z_1}\cdots\frac{dz_n}{z_n}
\\ && \makebox[4em]{}
=\prod_{j=1}^n\left(\frac{\eg(t^j;p,q)}{\eg(t;p,q)}
\prod_{1\leq m<s\leq 6}\eg(t^{j-1}t_mt_s;p,q)\right),
\label{eSelberg}\end{eqnarray}
where $|p|, |q|,|t|,|t_m| <1$ and $t^{2n-2}\prod_{m=1}^6t_m=pq$.
A conditional proof of this relation was given in \cite{die-spi:selberg}.
It depended on an evaluation of another elliptic hypergeometric integral
which was proven in \cite{rai:trans,spi:short}. In a special $p\to 0$
limit formula \eqref{eSelberg} reduces to a Gustafson multiple $q$-beta
integral \cite{gus:some} which, in turn, can be degenerated to the
Selberg integral (see, e.g., \cite{aar}). Integral \eqref{eSelberg} serves
as the measure for a very general class of biorthogonal functions
found by Rains \cite{rai:abelian,rai:trans} who defined a multivariable extension
of the author's two-index biorthogonal functions \cite{spi:theta2}
and elliptic analogues of the Koornwinder-Macdonald orthogonal polynomials
and their interpolation versions due to Okounkov.

There are many exact integration formulas or symmetry transformation
relations for elliptic hypergeometric integrals on roots systems
analogous to \eqref{eSelberg} or \eqref{E7}. A large
list of them can be found in \cite{SV2,SV3} where about a half
of the presented relations are formulated as conjectures.

\section{An elliptic Fourier transformation}

The Bailey chain technique discovered by Andrews \cite{andrews}
and Paule \cite{paule} is a tool for generating
infinite sequences  of identities for $q$-hypergeometric series
(see, e.g. \cite{aar}). It emerged from universalization of the proofs of
famous Rogers-Ramanujan identities, one of which has the following form
$$
\sum_{n=0}^\infty\frac{q^{n^2}}{(q;q)_n}=\frac{1}{(q;q^5)_\infty(q^4;q^5)_\infty}.
$$
The well-poised Bailey lemma \cite{and:bailey} provides the top known $q$-series
tool of such type. It served as an initial step for author's elliptic
generalization of this technique. First it was done for elliptic hypergeometric
series (see, e.g., \cite{spi:bailey1,war:extensions}). Subsequently it was generalized
to elliptic hypergeometric integrals \cite{spi:bailey}, which appeared
to be the very first extension of the Bailey chain technique to integrals.

Consider this construction, conditionally called an elliptic Fourier transformation,
in more detail. Take two functions $\alpha(z,t)$ and $\beta(w,t)$ depending
on a complex variable $z$ and a parameter $t\in\C$. They are said to form
an elliptic integral Bailey pair with respect to $t$ if they are
related by the following integral transformation
\begin{equation}
\beta(w,t)=M(t)_{wz}\alpha(z,t)
:=\frac{(p;p)_\infty(q;q)_\infty}{4\pi\textup{i}}\int_\mathbb{T}
\frac{\Gamma(tw^{\pm1}z^{\pm1};p,q)}{\Gamma(t^2,z^{\pm2};p,q)}
\alpha(z,t)\frac{dz}{z},
\label{baileypair}\end{equation}
where $|tw^{\pm1}|<1$. For wider region of parameters
this operator is defined by analytical continuation --- it
is necessary to replace the contour of integration $\T$
by a contour $C$ which separates the sequences of poles
converging to zero $z=tw^{\pm1}p^jq^k,\, j,k \in \Z_{\geq 0},$
from their reciprocals $z=t^{-1}w^{\pm1}p^{-j}q^{-k},
\, j,k \in \Z_{\geq 0},$ diverging to
infinity. Existence of such a contour is the only restriction
for the definition of analytically continued operator.

An integral analogue of the Bailey lemma has the
following form \cite{spi:bailey}.
\begin{theorem}
Let $\alpha(z,t)$ and $\beta(w,t)$ form an elliptic integral Bailey pair with respect
to $t$. Then the functions
$$
\alpha'(w,st)=D(s;u,w)\alpha(w,t),
$$
where
$$
D(s;u,w):=\Gamma(\sqrt{pq}s^{-1}u^{\pm1}w^{\pm1};p,q),
\qquad  D(s;u,w)D(s^{-1};u,w)=1,
$$
with $s, u\in\C$ being arbitrary new parameters, and
$$
\beta'(w,st)=D(t^{-1};u,w) M(s)_{wx}D(st;u,x)\beta(x,t)
$$
form an elliptic integral Bailey pair with respect to the parameter $st$.
\end{theorem}

After substitution of explicit expressions for $\alpha'$ and $\beta'$
into the required equality $\beta'(w,st)$ $=M(st)_{wz}\alpha'(z,st)$,
the proof boils down to the relation
\begin{equation}
M(s)_{wx}D(st;u,x)M(t)_{xz}=D(t;u,w)M(st)_{wz}D(s;u,z)
\label{operator_str}\end{equation}
which is equivalent to the elliptic beta integral evaluation formula \eqref{ellbeta}.
Relation \eqref{operator_str} is known as the operator form
of the star-triangle relation \cite{DM} formally depicted in the figure
(the black circle denotes integration over $x$ in the left-hand
side of \eqref{operator_str}).

\begin{figure}[ht]\vspace{0.3cm}
\begin{center}
\leavevmode \epsfxsize=5cm \epsffile{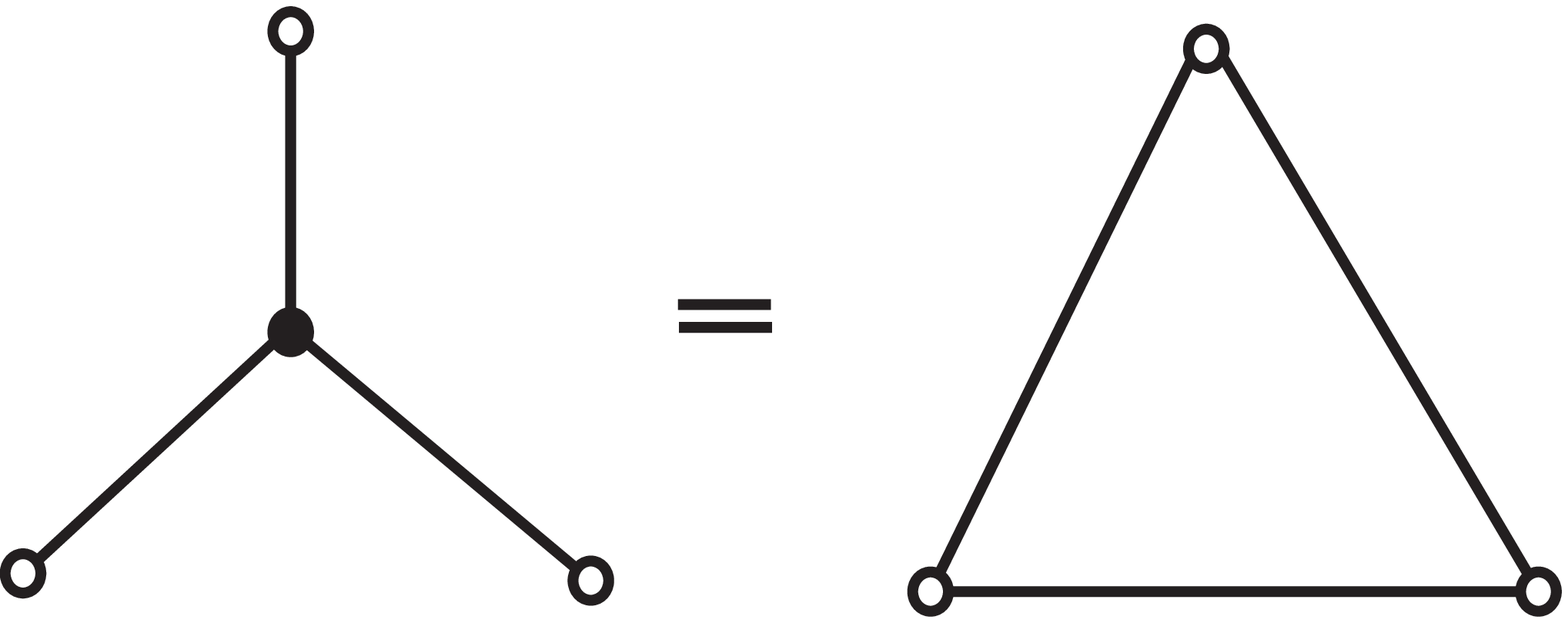}
\end{center}\vspace{-0.2cm}
\end{figure}

The main motivation for calling transformation \eqref{baileypair} an elliptic
analogue of the Fourier transformation comes from its inversion property
established in \cite{spi-war:inversions}. Namely, up to some
contour deformations the inversion relation is
equivalent to the reflection $t\to t^{-1}$:
$$
M(t)_{wz}M(t^{-1})_{zx}f(x)=f(w),
$$
or ``inversion = a sign change" (like in the Fourier transformation).
One can use theorem~2 for proving various identities for
elliptic hypergeometric integrals (e.g., of relation \eqref{E7}
or of the identity derived in \cite{bult}).

\section{The Yang-Baxter equation (YBE)}

The Yang-Baxter equation (YBE)  has the form \cite{bax:book,FT}
\begin{equation}
\mathbb{R}_{12} (u-v)\,\mathbb{R}_{13}(u)\, \mathbb{R}_{23}(v)
=\mathbb{R}_{23}(v)\,\mathbb{R}_{13}(u)\,\mathbb{R}_{12}(u-v),
\label{YBE}\end{equation}
where the operator $\mathbb{R}_{jk}(u)$ acts nontrivially
only in  $\V_j\otimes\V_k \subset
\V_1\otimes\V_2\otimes\V_3$ with $\V_j$ being some (in general different)
spaces. The variable $u\in\C$ is called the spectral parameter.

While investigating the eight-vertex model Baxter
found a YBE solution for $\V_j=\C^2$ (i.e., $\dim \V_j=2$)
\cite{bax}:
\begin{equation}
\mathbb{R}_{12}(u) = \sum_{a=0}^3 w_{a} (u)\,
\sigma_a \otimes\sigma_a,\quad w_{a}(u) = \frac{\theta_{a+1}
(u+\eta|\tau)}{\theta_{a+1}(\eta|\tau)},
\label{RBax}\end{equation}
where $\sigma_0=\II$ and $\sigma_{1,2,3}$ are the Pauli matrices,
$\theta_a(u|\tau)\equiv\theta_a(u)$ are the Jacobi theta functions.

Sklyanin \cite{skl} solved YBE for
$\dim \V_1=\dim \V_2=2, \dim \V_3=\infty,$ when YBE takes the form
$$
\mathbb{R}_{12} (u-v)\,\mathrm{L}_{13}(u)\, \mathrm{L}_{23}(v)
=\mathrm{L}_{23}(v)\,\mathrm{L}_{13}(u)\,\mathbb{R}_{12}(u-v)
$$
with $\mathbb{R}_{12} (u)$ being the Baxter R-matrix \eqref{RBax}.
In this case
\begin{eqnarray*} &&
\mathrm{L}_{13}(u):=\mathrm{L}(u):=
\sum_{a=0}^3 w_{a} (u)\, \sigma_a \otimes \mathbf{S}^a
\\ && \makebox[2em]{}
= \left(
\begin{array}{cc}
w_0(u)\,\mathbf{S}^0+w_3(u)\,\mathbf{S}^3 &
w_1(u)\,\mathbf{S}^1-\textup{i} w_2(u)\,\mathbf{S}^2 \\
w_1(u)\,\mathbf{S}^1+\textup{i} w_2(u)\,\mathbf{S}^2&
w_0(u)\,\mathbf{S}^0-w_3(u)\,\mathbf{S}^3
\end{array} \right),
\end{eqnarray*}
and the operators $\mathbf{S}^a$ generate the Sklyanin algebra:
 \begin{eqnarray}
\nonumber &&
\mathbf{S}^\alpha\,\mathbf{S}^\beta - \mathbf{S}^\beta\,\mathbf{S}^\alpha  =
\textup{i}\left(\mathbf{S}^0\,\mathbf{S}^\gamma +\mathbf{S}^\gamma\,
\mathbf{S}^0\right),
\\ &&
\mathbf{S}^0\,\mathbf{S}^\alpha - \mathbf{S}^\alpha\,\mathbf{S}^0  =
\textup{i}\mathbf{J}_{\beta \gamma}\left(\mathbf{S}^\beta\,\mathbf{S}^\gamma
+\mathbf{S}^\gamma\,\mathbf{S}^\beta\right),
\label{sklalg}
\end{eqnarray}
with $(\alpha,\beta,\gamma)$ a cycle of $(1,2,3)$ and the
structure constants $\mathbf{J}_{12}=\theta_1^2(\eta)\theta_4^2(\eta)/
\theta_2^2(\eta)\theta_3^2(\eta), $ etc. Relations \eqref{sklalg} define
an elliptic deformation of the $sl(2)$-algebra of rank 1.

The generators $\mathbf{S}^a$ can be realized explicitly as finite-difference operators
\cite{skl}
\begin{equation}
\left[\mathbf{S}^a\,\Phi\right](z) =\frac{\textup{i}^{\delta_{a,2}}
\theta_{a+1}(\eta)}{\theta_1(2 z) } \Bigl[\,\theta_{a+1} \left(2
z-2\eta\ell\right)\cdot \Phi(z+\eta)
- \theta_{a+1}
\left(-2z-2\eta\ell\right)\cdot \Phi(z-\eta)\, \Bigl],
\label{Sklgen}\end{equation}
where the variable $\ell\in\C$ is called the spin.

The relation between the Sklyanin
algebra and elliptic hypergeometric functions was first discussed
by Rains and Rosengren \cite{rai:abelian,ros:sklyanin}.
In \cite{S5} it was shown that the elliptic analogue of the
Euler-Gauss hypergeometric function can be derived as
a scalar product of solutions of generalized eigenvalue problems
for linear combinations of the Sklyanin algebra generators.

\section{A solution of the YBE for $\dim \V_j =\infty$}

The logical scheme of building an infinite-dimensional YBE solution of rank 1
by Derkachov and the author in \cite{DS}  consists of two steps:
\begin{itemize}
\item Take a defining $\mathrm{RLL}$-relation for $\dim \V_1=\dim \V_2=\infty,
\ \dim\V_3=2$ in the form
\begin{equation}
\mathbb{R}_{12}^?(u-v) \,\sigma_3 \mathrm{L}_1(u)
\,\sigma_3\,\mathrm{L}_2(v)=
\sigma_3\,\mathrm{L}_2(v) \,\sigma_3\mathrm{L}_1(u)\,\mathbb{R}_{12}^?(u-v),
\label{RLL}\end{equation}
where $\mathbb{R}_{12}^?(u)$ is an operator (the question mark means that it is unknown)
acting in the space of functions of two complex variables $\Phi(z_1,z_2)$,
$\mathrm{L}_1(u)=\mathrm{L}(u)$ with $z$ and $\ell$ replaced by $z_1$ and $\ell_1$,
$\mathrm{L}_2(v)=\mathrm{L}(v)$ with $z$ and $\ell$ replaced by $z_2$ and $\ell_2$.
Solve it using some auxiliary operators $\mathrm{S}_j, j=1,2,3$,
generating the permutation group $\mathfrak{S}_4$,
whose Coxeter relations are guaranteed by the elliptic beta integral.

\item Prove that the resulting $\mathbb{R}_{12}^?(u)$-operator
obeys the general YBE for $\dim\V_j=\infty$.
\end{itemize}

Extract from the R-matrix the permutation operator
$\mathbb{R}_{12}^?(u): = \mathbb{P}_{12}\,\mathrm{R}_{12}(u)$.
Then relation \eqref{RLL} takes the form
\begin{equation}
\mathrm{R}_{12}(u-v)\,\mathrm{L}_1(u_1,u_2)\,\sigma_3\,\mathrm{L}_2(v_1,v_2)=
\mathrm{L}_1(v_1,v_2)\,\sigma_3\,\mathrm{L}_2(u_1,u_2)\,\mathrm{R}_{12}(u-v),
\label{RLL'}\end{equation}
where
$$
u_1 =
\frac{u}{2}+ \eta\,(\ell_1+\frac{1}{2}),\
u_2 = \frac{u}{2}-\eta\,(\ell_1+\frac{1}{2}),\
v_1 = \frac{v}{2} + \eta\,(\ell_2+\frac{1}{2}),\
v_2 = \frac{v}{2}-\eta\,(\ell_2+\frac{1}{2}).
$$
Notice that in \eqref{RLL'} the operator
$\mathrm{R}_{12}(u-v)\equiv \mathrm{R}_{12}(u_1,u_2|v_1,v_2)$
just permutes parameters in the product of $\mathrm{L}$-operators.
Denote
$$
\mathbf{u}\equiv(u_1,u_2,v_1,v_2),
\quad s\mathbf{u} =(v_1,v_2,u_1,u_2), \qquad s = s_2 s_1s_3s_2,
$$
where $s_i$ are elementary permutations
$$
s_{1}\mathbf{u} = (u_2,u_1,v_1,v_2), \quad
s_{2}\mathbf{u}  = (u_1,v_1,u_1,v_2), \quad
s_{3}\mathbf{u} = (u_1,u_2,v_2,v_1).
$$
Define $\mathrm{S}_j$-operators by the relations
\begin{eqnarray*} &&
\mathrm{S}_1(\mathbf{u})\,\mathrm{L}_1(u_1,u_2) =
\mathrm{L}_1(u_2,u_1)\,\mathrm{S}_1(\mathbf{u}),\quad
\mathrm{S}_3(\mathbf{u})\,\mathrm{L}_2(v_1,v_2) =
\mathrm{L}_2(v_2,v_1)\,\mathrm{S}_3(\mathbf{u}),
\\ &&
\mathrm{S}_2(\mathbf{u})\,\mathrm{L}_1(u_1,u_2)\,\sigma_3\,\mathrm{L}_2(v_1,v_2)=
\mathrm{L}_1(u_1,v_1)\,\sigma_3\,\mathrm{L}_2(u_2,v_2)\,\mathrm{S}_2(\mathbf{u}).
\end{eqnarray*}
These operators generate the permutation group $\mathfrak{S}_4$ if
they satisfy the Coxeter relations:
\begin{equation}
\mathrm{S}_i^2=\II, \ \
\mathrm{S}_i\mathrm{S}_j=\mathrm{S}_j\mathrm{S}_i, \ |i-j|>1, \ \
\mathrm{S}_j\mathrm{S}_{j+1}\mathrm{S}_j
=\mathrm{S}_{j+1}\mathrm{S}_j\mathrm{S}_{j+1},
\label{Coxeter}\end{equation}
where we assume the following multiplication rule
$\mathrm{S}_j\mathrm{S}_k:=\mathrm{S}_j(s_k\mathbf{u})
\mathrm{S}_k(\mathbf{u})$.

\begin{theorem} The operator
$$
\mathrm{R}_{12}(\mathbf{u}) =
\mathrm{S}_2(s_1s_3s_2\mathbf{u})\,
\mathrm{S}_1(s_3s_2\mathbf{u})\,\mathrm{S}_3(s_2\mathbf{u})\,
\mathrm{S}_2(\mathbf{u})
$$
solves the initial $\mathrm{RLL=LLR}$ relation \eqref{RLL'}.
\end{theorem}

The proof is straightforward and follows from the intertwining
properties of $\mathrm{S}_j$-operators.

\begin{theorem} The operator $\mathbb{R}_{12}(u)=\mathbb{P}_{12}\,\mathrm{R}_{12}(\mathbf{u})$
solves the YBE \eqref{YBE}.
\end{theorem}

The following permutation of parameters in the product
of three $\mathrm{L}$-operators
$$
\mathrm{L}_1 (u_1,u_2)\,\sigma_3\, \mathrm{L}_2(v_1,v_2)\,\sigma_3\,
\mathrm{L}_3(w_1,w_2)
\to \mathrm{L}_1 (w_1,w_2)\,\sigma_3\,
\mathrm{L}_2(v_1,v_2) \,\sigma_3\,\mathrm{L}_3(u_1,u_2)
$$
can be realized in two different ways indicating that
\begin{eqnarray}\nonumber &&
\mathrm{R}_{12}(v_{1},v_2|w_{1},w_2)\,\mathrm{R}_{23}(u_{1},u_2|w_1,w_2)\,
\mathrm{R}_{12}(u_1,u_2|v_{1},v_2)
\\ &&\makebox[2em]{}
=\mathrm{R}_{23}(u_{1},u_{2}|v_1,v_2)\,\mathrm{R}_{12}(u_1,u_{2}|w_1,w_{2})\,
\mathrm{R}_{23}(v_{1},v_2|w_{1},w_2).
\label{YBErel}\end{eqnarray}
This relation is proved directly with the help of cubic Coxeter
relations for $\mathrm{S}_j$-operators alone, i.e. it is just a word
identity in the group algebra of the braid group $\mathfrak{B}_4$
(i.e., $\mathrm{S}_j^2=\II$ relations are not used).
Multiplying relation \eqref{YBErel} by the appropriate
permutation operators $\mathbb{P}_{ij}$ one comes to the original
YBE \eqref{YBE}.

As to the explicit construction of needed operators $\mathrm{S}_j$,
a miracle takes place --- after demanding that $\mathrm{S}_{1,2}$
are integral operators one comes to the elliptic Fourier
transformation described above!

Denote $p=e^{2\pi\textup{i}\tau}$ and $q=e^{4\pi\textup{i}\eta}$.
Then, for $|p|, |q|<1$ and a special choice of periodic factors emerging from
solutions of finite difference equations one obtains
\begin{eqnarray*} &&
[\mathrm{S}_2\Phi](z_1,z_2)=
D(e^{2\pi\textup{i}(v_1-u_2)};e^{2\pi\textup{i}z_1},
e^{2\pi\textup{i}z_2})\Phi(z_1,z_2),
\\ &&
[\mathrm{S}_1\Phi](z_1,z_2) = e^{-\pi \textup{i} z_1^2/\eta}
M(e^{2\pi \textup{i}(u_2-u_1)})_{e^{2\pi\textup{i}z_1},e^{2\pi\textup{i}z}}
e^{\pi \textup{i} z^2/\eta}\Phi(z,z_2),
\end{eqnarray*}
where $D(t;x,y)$ and $M(t)_{xy}$ are the elliptic integral
Bailey lemma entries. Notice the reduced form of the parameter dependence
$\mathrm{S}_2({\bf u})=\mathrm{S}_2(u_2-v_1)$ and
$\mathrm{S}_1({\bf u})=\mathrm{S}_1(u_1-u_2)$.
The operator $\mathrm{S}_3({\bf u})$ has the same form as
$\mathrm{S}_1({\bf u})$ with $z_1$ replaced by $z_2$ and $u_1-u_2$ by $v_1-v_2$.
In this picture the Coxeter relations
coincide with the identities for $D$ and $M$-operators.
In particular, the cubic Coxeter relation $\mathrm{S}_1\mathrm{S}_2\mathrm{S}_1=
\mathrm{S}_2\mathrm{S}_1\mathrm{S}_2$ is guaranteed by the
elliptic beta integral since it has the explicit form
$$
\mathrm{S}_1(a)\,\mathrm{S}_2(a+b)\,\mathrm{S}_1(b)=
\mathrm{S}_2(b)\,\mathrm{S}_1(a+b)\,\mathrm{S}_2(a),
$$
coinciding with equality \eqref{operator_str}.
The relation $\mathrm{S}_1^2=\II$ explicitly looks as
$\mathrm{S}_1(-a)\mathrm{S}_1(a)=\II$, which is the inversion
relation for the elliptic Fourier transformation. Thus,
$$
\text{elliptic beta integral} = \text{star-triangle relation}
= \text{cubic Coxeter relation.}
$$

After the similarity transformation removing the exponentials
$e^{\pi \textup{i} z^2/\eta}$, which is
equivalent to passing to the Sklyanin algebra generators
$\mathbf{S}^a\to e^{\pi \textup{i} z^2/\eta}\mathbf{S}^a
e^{-\pi \textup{i} z^2/\eta}$, one obtains the general
R-operator as a double integral operator with an elliptic
hypergeometric kernel of the form
\begin{eqnarray}  \nonumber &&
[\mathbb{R}_{12}(\mathbf{u})f](x_1,x_2)
=\frac{(p;p)_\infty^2(q;q)_\infty^2}{(4\pi \textup{i})^2}
\Gamma(\sqrt{pq}x_1^{\pm1}x_2^{\pm 1}e^{2\pi\textup{i}(v_2-u_1)};p,q)
\\  \nonumber && \makebox[2em]{} \times
\int_{\mathbb{T}^2}
\frac{
\Gamma(e^{2\pi\textup{i}(v_1-u_1)}x_2^{\pm1}x^{\pm1},
e^{2\pi\textup{i}(v_2-u_2)}x_1^{\pm1}y^{\pm1};p,q)}
{\Gamma(e^{4\pi\textup{i}(v_1-u_1)},e^{4\pi\textup{i}(v_2-u_2)},
x^{\pm2},y^{\pm2};p,q)}
\\ && \makebox[6em]{} \times
\Gamma(\sqrt{pq}e^{2\pi\textup{i}(v_1-u_2)}x^{\pm1}y^{\pm1};p,q)
f(x,y)\frac{dx}{x}\frac{dy}{y},
\label{YBEsol}\end{eqnarray}
where one should impose some mild constraints on the parameter values
in order to satisfy the intertwining relations with  the taken integration
contour $\mathbb T$ \cite{DS}.

\section{The elliptic modular double}

The derived R-operator \eqref{YBEsol} is symmetric in $p$ and $q$. Hence there exists
a second RLL-relation obtained from the first one \eqref{RLL} by permutation of
$p$ and $q$:
$$
\mathbb{R}_{12}(u-v)\,\sigma_3 \mathrm{L}_1'(u)
\,\sigma_3\,\mathrm{L}_2'(v)=
\sigma_3\,\mathrm{L}_2'(v) \,\sigma_3\mathrm{L}_1'(u)\,\mathbb{R}_{12}(u-v),
$$
where
$\mathrm{L}'(u)=\mathrm{L}(\text{fixed } u, \text{ fixed } g=\eta(2\ell+1), \
2\eta\leftrightarrow \tau).$
This means that there is a second copy of the Sklyanin algebra generated
by the operators
\begin{eqnarray*}&&
\mathbf{\tilde S}^a =
e^{\frac{2\pi i}{\tau} z^2 }\, \frac{\textup{i}^{\delta_{a,2}}
\theta_{a+1}(\frac{\tau}{2})|2\eta)}{\theta_1(2 z|2\eta) }
\Bigl[\,\theta_{a+1} \left(2z-g+\frac{\tau}{2}|2\eta\right)
e^{\frac{1}{2}\tau \partial_z}
\\ && \makebox[2em]{}
 - \theta_{a+1} \left(-2z-g+\frac{\tau}{2}|2\eta\right)
e^{-\frac{1}{2}\tau \partial_z}\, \Bigl]
e^{-\frac{2\pi i}{\tau} z^2}, \quad
e^{\alpha\partial_z}f(z)=f(z+\alpha).
\end{eqnarray*}

The direct product of two such Sklyanin algebras
was introduced by the author in \cite{S5} under the name
``elliptic modular double". It represents an elliptic generalization
of the Faddeev modular double for the $sl_q(2)$ quantum algebra \cite{fad:mod}.
Vice versa, demanding existence of  the elliptic modular double
together with the meromorphy of functions  in the variable $e^{2\pi\textup{i}z}$
removes periodic factors in solutions of difference equations
and determines the operators $\mathrm{S}_j$ uniquely.

\section{The superconformal index}

The first physical interpretation of elliptic hypergeometric integrals
has been found by the author in the context of
Calogero-Sutherland type models \cite{spi:cs}.
The most remarkable application of such integrals in physics
has been discovered by Dolan and Osborn \cite{DO}.

The superconformal index \cite{Kinney,rom} is a topological index of four-dimensional ($4d$)
supersymmetric gauge field theories with local gauge invariance group
$G$ and global flavor symmetry group $F$
and some set of fields described by irreducible representations
of these groups $R_{G, j}$ and $R_{F,j}$. A heuristic
derivation of this object resulted in the following matrix integral
\begin{equation}
I(y;p,q) = \int_{G} d \mu(z)\,
\exp \bigg ( \sum_{n=1}^{\infty}
\frac 1n \text{ind}\big(p^n ,q^n, z^n , y^ n\big ) \bigg )
\label{SCI}\end{equation}
with Haar measure $d\mu(z)$ and
\begin{eqnarray}\nonumber &&
\text{ind}(p,q,z,y) =  \frac{2pq - p - q}{(1-p)(1-q)} \chi_{adj_G}(z)\cr
\\ && \makebox[2em]{}
+ \sum_j \frac{(pq)^{R_j/2}\chi_{r_F,j}(y)\chi_{r_G,j}(z) - (pq)^{1-R_j/2}
\chi_{{\bar r}_F,j}(y)\chi_{{\bar r}_G,j}(z)}{(1-p)(1-q)}.
\label{ind}\end{eqnarray}
Here $\chi_{R_F,j}(y)$ and $\chi_{R_G,j}(z)$ are characters
of the respective representations, and $R_j$ are some fractional numbers
(R-charges of the fields).
For instance, for the unitary group $SU(N)$,
$z=(z_1,\ldots,z_N), \ \prod_{a=1}^Nz_a=1$, one has
$$
\int_{SU(N)} d\mu(z) \ = \   \frac{1}{N!} \int_{\mathbb{T}^{N-1}}
\Delta(z) \Delta(z^{-1}) \prod_{a=1}^{N-1} \frac{dz_a}{2 \pi \textup{i} z_a},
$$
where $\Delta(z) \ = \ \prod_{1 \leq a < b \leq N} (z_a-z_b)$
is the Vandermonde determinant.

Where is the elliptic beta integral here?
Let us take $G=SU(2),$ $F=SU(6)$, and the representations
(``$adj$"=adjoint, ``$f$"=fundamental)
\begin{eqnarray*} &&
1)\ \text{``vector superfield": } (adj, 1),
\quad \chi_{SU(2),adj}(z)=z^2+z^{-2}+1,
 \\ &&
2)\ \text{``chiral superfield": } (f, f),
\quad \chi_{SU(2),f}(z)=z+z^{-1},\qquad R_f=1/3,
\\ && \makebox[0em]{}
\chi_{SU(6),f}(y)=\sum_{k=1}^6y_k,
\quad \chi_{SU(6),\bar f}(y)=\sum_{k=1}^6y_k^{-1},
\quad \prod_{k=1}^6y_k=1,
\end{eqnarray*}
 Then after denoting $t_k=(pq)^{1/6}y_k$, $k=1,\ldots,6,$ the superconformal
index formula reproduces identically the left-hand side
of the elliptic beta integral \eqref{ellbeta}  \cite{DO}.

In order to generate the right-hand side expression
in \eqref{ellbeta} one should take $G=1,$ $F=SU(6)$ with the single
``chiral superfield" described by antisymmetric tensor of the second rank
$T_A:\quad \Phi_{ij}=-\Phi_{ji},$
$$
\chi_{SU(6),T_A}(y)=\sum_{1\leq i<j\leq 6}y_iy_j,
\qquad R_{T_A}=2/3.
$$
So, the elliptic beta integral evaluation formula
shows that two functions on characters \eqref{SCI}
for different sets of representations of two different $G\times F$-groups
coincide. This equality of superconformal indices describes the confinement
phenomenon in the simplest supersymmetric quantum chromodynamics model
or the Seiberg duality \cite{seiberg}. An almost complete list of
such dualities for simple gauge groups and corresponding elliptic
hypergeometric integral identities (about half of which are conjectures)
emerging as equalities of dual indices is given in \cite{SV2,SV3}.
In general, elliptic hypergeometric integrals define new matrix models
and describe the most complicated known class of computable
nonperturbative path integrals in four-dimensional quantum field theories.

\section{A $4d/2d$  correspondence}

Solutions of the YBE are related to integrable spin chain models
and $2d$ Ising type spin systems \cite{FT}.
Let us replace in the Coxeter/Bailey relation \eqref{operator_str}
$u,x,z \to e^{\textup{i}
u},  e^{\textup{i}x},  e^{\textup{i}z}$ with real $u,x,z$
and act by it on a localized
``continuous spin" state function $(\delta(z-y)+\delta(z+y))/2$,
where $\delta(z)$ is the Dirac delta function.
This yields the elliptic beta integral rewritten in the form
\begin{eqnarray*}
&& \int_{0}^{2\pi} \rho(u)D_{\xi-\alpha}(x,u)D_{\alpha+\gamma}(y,u)
D_{\xi-\gamma}(w,u)du
=\chi D_{\alpha}(y,w)D_{\xi-\alpha- \gamma}(x,w)D_\gamma(x,y),
\end{eqnarray*}
where
\begin{eqnarray*}
&& D_\alpha(x,y)=D(e^{-\alpha}; e^{\textup{i} x},e^{\textup{i} y}),
\quad \rho(u)=\frac{(p;p)_\infty (q;q)_\infty}{4\pi}
\theta(e^{2 \textup{i}  u};p)\theta(e^{-2 \textup{i}  u};q),
\\  \makebox[0em]{}
&& \chi=  \Gamma(e^{-2\alpha},e^{-2\gamma},e^{2\alpha+2\gamma-2\xi};p,q),
\qquad e^{-\xi}=\sqrt{pq},
\end{eqnarray*}
which coincides with the functional star-triangle relation considered
by Bazhanov and Sergeev \cite{BS}.

Now one can interpret the figure given earlier as a transformation
of the elementary cell partition function: circles
carry spins $u,w,\ldots$ with the self-energy $\rho(u)$,
edges carry Boltzmann weights $D_\alpha$, the
black circle contains the integration (summation) over $u$-spin values.
Using such a relation one can map the honeycomb, triangular, and
square lattice spin models onto each other and
their partition functions become equal to some (multiple) elliptic
hypergeometric integrals.

The same figure can be interpreted in the context of quiver gauge
theories: black circles denote the gauge group $G$
(in superconformal indices this corresponds to the contribution of
vector superfields in the adjoint representation and integration over
the gauge group),  white circles
denote some external flavor groups, the edges denote the bifundamental
fields (their contributions to superconformal indices yield
the Boltzmann weights $D_\alpha$), etc. The presence of more than one
black circle would indicate that the gauge group is not simple.
Evidently one can construct in this way lattice-type $4d$ quiver
field theories so that their superconformal indices coincide with
partition functions of various statistical mechanics models
with the continuous spin values \cite{stat} since both are described by
the same elliptic hypergeometric integrals.

As shown in \cite{stat}, one can go further and interpret
the symmetry transformations for elliptic hypergeometric
integrals on root systems as the star-star relations \cite{bax_ss}.
Conjecturally, to each such non-trivial relation one can associate
an elliptic integrable system.
As a consequence, the Seiberg duality of $4d$ field theories
becomes related to the  Kramers-Wannier type duality
transformations for partition functions of $2d$ spin systems.

\section{Conclusion}

We can conclude that the elliptic hypergeometric functions are
universal objects with wide applications. This brief survey
does not cover all their known instances and, in particular,
the list of given references is incomplete (for its extensions,
see \cite{spi:umnrev,SV2,SV3}). We finish by
indicating various fields where elliptic hypergeometric
functions have found their applications and more is expected to lie ahead.

{\bf In mathematics:}
theory of analytic finite-difference equations (e.g., the
elliptic hypergeometric equation), harmonic analysis on root
systems, representation theory, theory of $SL(3,\Z)$ automorphic forms,
approximation theory, continued fractions, combinatorics, topology, etc.

{\bf In theoretical and mathematical physics:}
$4d$ supersymmetric dualities,
integrable $N$-particle quantum mechanical systems,
$2d$ topological field theories \cite{gprr},
$2d$ solvable models of statistical mechanics and
noncomact spin chains, random matrices and
determinantal point processes \cite{bgr}, etc.

\medskip

This paper is authentic with author's talk given at the conference
dedicated to 125th anniversary of S. Ramanujan in Delhi.
The author is deeply indebted to organizers
for the invitation to this meeting and kind hospitality during it.
This work is supported in part by RFBR grant no.  12-01-00242.

\end{document}